\documentclass[english]{amsart}

\usepackage{babel}
\usepackage{amstext}
\usepackage{amsmath}
\usepackage{amsfonts}
\usepackage{latexsym}
\usepackage{ifthen}
\usepackage[all,2cell,dvips]{xy} \UseAllTwocells \SilentMatrices
\xyoption{all}
\newcommand{\codim}{{\rm codim}}

\newtheorem{lemma1}[equation]{}

\newenvironment{lemma}{\begin{lemma1}{\bf Lemma.}}{\end{lemma1}}

\newenvironment{theorem}{\begin{lemma1}{\bf Theorem.}}{\end{lemma1}}

\newenvironment{definition}{\begin{lemma1}{\bf Definition.}}{\end{lemma1}}

\newenvironment{conjecture}{\begin {lemma1}{\bf Conjecture.}}{\end{lemma1}}

\newcommand{\Z}{\ensuremath{\mathbb{Z}}}
\newcommand{\C}{\ensuremath{\mathbb{C}}}

\newcommand{\PP}{\ensuremath{\mathbb{P}}}

\newcommand{\Alb}[1]{\ensuremath{\mbox{Alb}(#1)\hspace{0.5ex}}}

\newcommand{\holom}[3]{\ensuremath{#1:#2  \rightarrow #3}}

\makeatletter
\ifnum\@ptsize=0  \fi \ifnum\@ptsize=2  \fi 

\newcommand\sF{{\mathcal F}}

\newcommand\sO{{\mathcal O}}

\DeclareMathOperator*{\Pic0}{Pic^0}
\DeclareMathOperator*{\Nm}{Nm}

\DeclareMathOperator*{\num}{num}

\title {$M$-regularity of the Fano surface} 
\date{4th April, 2007}

\author{Andreas H\"oring}
\address{Andreas H\"oring, IRMA, Universit\'e Louis Pasteur, 7 rue Ren\'e Descartes, 67084 Strasbourg, France}
\email{andreas.hoering@ujf-grenoble.fr}

\begin{document}

\begin{abstract}
In this note we show that the Fano surface in the intermediate Jacobian of a smooth 
cubic threefold is $M$-regular in the sense of Pareschi and Popa. 
\end{abstract}

\maketitle


\section{Introduction}

Let $X^3 \subset \PP^4$ be a smooth cubic threefold, then its intermediate Jacobian
\[
J(X) := H^{2,1}(X, \C)^* / H_3(X, \Z)
\]
is a principally polarised abelian variety $(J(X), \Theta)$
of dimension five that is not a Jacobian of a curve
\cite[Thm.0.12]{CG72}.
The Fano scheme $F$ parametrising lines contained in $X$ is a smooth surface, and the
Abel-Jacobi map $F \rightarrow J(X)$ is an embedding 
that induces an isomorphism $\Alb{F)} \simeq J(X)$
\cite[Thm.0.6,0.9]{CG72}. Furthermore the cohomology class of $F \subset J(X)$ is minimal,
that is 
\[
[F] = \frac{\Theta^3}{3!}.
\]
There is only one other known family of 
examples of principally polarised abelian varieties $(A, \Theta)$ of dimension $n$ such that
for $1 \leq d \leq n-2$, a minimal cohomology class $\frac{\Theta^{n-d}}{(n-d)!}$ can be 
represented by an effective cycle
of dimension $d$: the Jacobians of curves $J(C)$ where the suvarieties 
$W_d(C) \subset J(C)$ have minmal cohomology class.
O. Debarre has shown that on a Jacobian these are the only subvarieties having minimal class
\cite[Thm.5.1]{Deb95}, 
furthermore by a theorem of Z. Ran \cite[Thm.5]{Ran81}, the only principally polarised abelian fourfolds 
with a subvariety of minimal class are (products of) Jacobians of curves. 
In higher dimension few things are known about subvarieties having minimal class.

In \cite{PP06}, G. Pareschi and M. Popa introduce a new approach to the characterisation 
of these subvarieties: they consider the  
(probably more tractable) cohomological properties of the twisted structure sheaf of the subvariety.
More precisely we have the following conjecture.

\begin{conjecture} \cite{Deb95},\cite{PP06}  \label{conjecturePP}
Let $(A, \Theta)$ be an irreducible principally polarised abelian varieties of dimension $n$, 
and let $Y$ be a nondegenerate subvariety (cf. \cite[p.464]{Ran81}) of $A$ of dimension $d\le n-2$.
The following statements are equivalent.
\begin{enumerate}
\item The variety $Y$ has minimal cohomology class, i.e. $[Y] = \frac{\Theta^{n-d}}{(n-d)!}$.
\item The twisted structure sheaf $\sO_Y(\Theta)$ is $M$-regular (cf. definition \ref{definitiongvsheaf} below), 
and $h^0 (Y, \sO_Y(\Theta) \otimes P_\xi) = 1$ for 
$P_\xi \in \Pic0(A)$ general.
\item Either $(A, \Theta)$ is the Jacobian of a curve of genus $n$ and
$Y$ is a translate of $W_d(C)$ or $-W_d(C)$, or $n= 5$, $d=2$ and $(A, \Theta)$ is the
intermediate Jacobian of a smooth cubic threefold and $Y$ is a translate of $F$ or $-F$. 
\end{enumerate}
\end{conjecture} 

The implication $2) \Rightarrow 1)$ is the object of \cite[Thm.B]{PP06}.
The implication $3) \Rightarrow 2)$ has been shown for Jacobians of curves in \cite[Prop.4.4]{PP03}.
We complete the proof of this implication by treating the case of the intermediate Jacobian. 

\begin{theorem}
\label{theoremfanosurfacegvsheaf}
Let $X^3 \subset \PP^4$ be a smooth cubic threefold, and let $(J(X), \Theta)$ be its intermediate Jacobian. 
Let $F \subset J(X)$ be an Abel-Jacobi embedded copy 
of the Fano variety of lines in $X$. Then $\sO_F(\Theta)$ is $M$-regular and
$h^0(F, \sO_F(\Theta) \otimes P_\xi)=1$ for $P_\xi \in \Pic0 J(X)$ general. 
\end{theorem}

Since the properties considered are invariant under isomorphisms, the theorem implies the same statement
for $-F$.

The study of the remaining open implications of conjecture \ref{conjecturePP} is a much harder task
than the proof of theorem \ref{theoremfanosurfacegvsheaf}. 
In an upcoming paper we will start to investigate this problem under the
additional hypothesis that $(A, \Theta)$ is the intermediate Jacobian of a generic smooth cubic threefold. 
In this case we can show the following statement. 

\begin{theorem} \cite{H07}
Let $X^3 \subset \PP^4$ be a {\em general} smooth cubic threefold. Let $(J(X), \Theta)$ be its intermediate Jacobian, and let $F \subset J(X)$ be an Abel-Jacobi embedded copy 
of the Fano variety of lines in $X$. Let $S \subset J(X)$ be a surface that has minimal cohomology class,
i.e. $[S] = \frac{\Theta^{3}}{3!}$.
Then $S$ is a translate of $F$ or $-F$.
\end{theorem}

\begin{center}
{\bf Notation and basic facts.}
\end{center}

We work over an algebraically closed field of characteristic different from 2.
We will denote by $D \equiv D'$ the linear equivalence of divisors, and by
$D \equiv_{\num} D'$ the numerical equivalence.

For $(A, \Theta)$ a principally polarised abelian variety (ppav), 
we identify $A$ with $\hat{A}=\Pic0(A)$ via the morphism induced by $\Theta$.
If $\xi \in A$ is a point, we denote by $P_\xi$ the corresponding point in $\hat{A}=\Pic0(A)$
which we consider as a numerically trivial line bundle on $A$.

\begin{definition} \cite{PP06b}
\label{definitiongvsheaf}
Let $(A, \Theta)$ be a ppav 
of dimension $n$, and let $\sF$ be a coherent sheaf
on $A$. For all $n \geq i > 0$, we denote by
\[
V^i_\sF := \{ \xi \in A \ | \ h^i(A, \sF \otimes P_\xi)>0 \}
\]
the $i$-th cohomological support locus of $\sF$. 
We say that $\sF$ is $M$-regular if
\[
\codim V^i_\sF > i 
\]
for all $i \in \{ 1, \ldots, n \}$.
\end{definition}

If $l \subset X$ is a line, we will denote by $[l]$ the corresponding point of the Fano surface $F$
and by $D_l \subset F$ the incidence curve of $l$,
that is, $D_l$ parametrises lines in $X$ that meet $l$. Furthermore we have by
\cite[\S 10]{CG72}, \cite[\S 6]{We81} and Riemann-Roch that
\begin{eqnarray}
\label{Theta2C}
\sO_F(\Theta) & \equiv_{\num} & 2 D_l,
\\
\label{KF3C}
K_F & \equiv_{\num} & 3 D_l,
\\
\label{Dintersection}
D_l \cdot D_l &=& 5,
\\
\label{chitheta}
\chi(F, \sO_F(\Theta)) &=& 1. 
\end{eqnarray}

\section{Prym construction of the Fano surface}
\label{sectionprym}

We recall the construction of the Fano surface as a special subvariety of a Prym variety \cite{Bea82, Bea82b}:
let $\tilde{C}:=D_{l_0} \subset F$ be the incidence curve of a general line $l_0 \subset X$.
Let $X'$ be the blow-up of $X$ in $l_0$. Then the projection from $l_0$ induces
a conic bundle structure $X' \rightarrow \PP^2$
with branch locus $C \subset \PP^2$ a smooth quintic. This conic bundle induces 
a natural connected
\'etale covering of degree two $\holom{\pi}{\tilde{C}}{C}$ (cf. \cite[Ch.I]{Bea77} for details),
and we denote by \holom{\sigma}{\tilde{C}}{\tilde{C}}
the involution induced by $\pi$.

The kernel of the 
norm morphism \holom{\Nm}{J \tilde{C}}{JC} has two connected components which we will denote by
$P$ and $P_1$. The zero component $P$ is called the Prym variety associated to $\pi$, and it is
isomorphic as a ppav to $J(X)$ \cite[Thm.2.1]{Bea77}. 

Let $H \subset C$ be an effective 
divisor given by a hyperplane section in $\PP^2$. Then $H$ has degree five
and $h^0(C, \sO_C(H))=3$, so the complete linear system $g_5^2$ corresponds to a $\PP^2 \subset C^{(5)}$.  
We choose a divisor $\tilde{H} \in \tilde{C}^{(5)}$ such that $\pi^{(5)}([\tilde{H}]) = [H]$, where
\holom{\pi^{(5)}}{\tilde{C}^{(5)}}{C^{(5)}} is the morphism induced by $\pi$ on the symmetric products.
Let \holom{\phi_H}{C^{(5)}}{J C} and \holom{\phi_{\tilde{H}}}{\tilde{C}^{(5)}}{J \tilde{C}} 
be the Abel-Jacobi maps given by $H$ and $\tilde{H}$. We have a commutative diagram
\[
\xymatrix{
\tilde{C}^{(5)} \ar[r]^{\phi_{\tilde{H}}} \ar[d]_{\pi^{(5)}}
& J \tilde{C} \ar[d]^{\Nm}
\\
C^{(5)} \ar[r]_{\phi_H} & J C
}
\]
The fibre of $\phi_{\tilde{H}} (\tilde{C}^{(5)}) \rightarrow \phi_H (C^{(5)})$ over the point $0$
(and thus the intersection of  $\phi_{\tilde{H}} (\tilde{C}^{(5)})$ with $\ker \Nm$) has two connected
components $F_0 \subset P$ and $F_1 \subset P_1$.
If we identify $P$ and $P_1$ via $\tilde{H} - \sigma(\tilde{H})$, 
we obtain an identification $F_1=-F_0$ \cite[p.360]{Bea82}.
The (non-canonical) isomorphism of ppavs
$P \simeq J(X)$ transforms $F_0$ into a translate of the Fano surface $F$
\cite[Thm.4]{Bea82}.

{\it From now on we will identify $P$ (resp. $F_0$) and $J(X)$ (resp. some Abel-Jacobi emdedded copy
of the Fano surface $F$)}.

We will now prove two technical lemmata on certain linear systems on $\tilde{C}$.
The first is merely a reformulation of  \cite[\S 2,ii)]{Bea82b}.

\begin{lemma}
The line bundle $\sO_{\tilde{C}}(\tilde{C})$ is a base-point free pencil of degree five such that 
any divisor $D \in |\sO_{\tilde{C}}(\tilde{C})|$ satisfies $\pi_* D \equiv H$. 
\end{lemma}

{\bf Proof.} 
We define a morphism \holom{\mu}{\tilde{C}=D_{l_0}}{l_0\simeq \PP^1}
by sending $[l] \in \tilde{C}$ to $l \cap l_0$. Since $l_0$ is general and through
a general point of $l_0$ there are five lines distinct from $l_0$, 
the morphism $\mu$ has degree 5. 
If $[l] \in F$, then $D_l \cdot D_{l_0}=5$ by formula (\ref{Dintersection}), so for $[l] \neq [l_0]$ 
the divisor $D_{l_0} \cap D_l \in |\sO_{\tilde{C}}(D_l)|$
is effective. Furthermore $\pi_* D_{l} \equiv H$, since $\pi_* D_{l}$ is the intersection of $C \subset \PP^2$
with the image of $l$ under the projection $X' \rightarrow \PP^2$. 
By specialisation the linear system $|\sO_{\tilde{C}}(\tilde{C})|$ is not empty
and a general divisor $D$ in it 
corresponds to the five lines distinct from $l_0$
passing through a general point of $l_0$. Hence
$\sO_{\tilde{C}}(\tilde{C}) \simeq \mu^* \sO_{\PP^1}(1)$ and $\pi_* D \equiv H$. $\square$ 

\begin{lemma}
\label{lemmadegreefiveten}
The sets 
\[
V_0' := \{
\xi \in P \ | \ h^0(\tilde{C}, \sO_{\tilde{C}}(\tilde{C}) \otimes P_\xi) > 0
\}
\] 
\[
V_1' := \{
\xi \in P \ | \ h^0(\tilde{C}, \sO_{\tilde{C}}(2 \tilde{C}) \otimes P_\xi) > 1
\}
\] 
are contained in translates of $F \cup (-F)$.
\end{lemma}

{\bf Proof.} 

1) Let $D \in |\sO_{\tilde{C}}(\tilde{C}) \otimes P_\xi|$ 
be an effective divisor. Then $\pi_* \tilde{C} \equiv \pi_* D \equiv H$. It follows that 
$D \in (\phi_{\tilde{H}} (\tilde{C}^{(5)}) \cap \ker \Nm)$, so $D$ is in $F$ or $-F$. 

2) We follow the argument in \cite[\S 3]{Bea82b}.
By \cite[\S 2,iv)]{Bea82b} we have $h^0(\tilde{C}, \sO_{\tilde{C}}(\tilde{C}+\sigma(\tilde{C})))=4$, 
so $h^0(\tilde{C}, \sO_{\tilde{C}}(2 \tilde{C}))$ is odd. It follows from the deformation invariance of
the parity \cite[p.186f]{Mum71} that
\[
V_1' = \{
\xi \in P \ | \ h^0(\tilde{C}, \sO_{\tilde{C}}(2 \tilde{C}) \otimes P_\xi) \geq 3
\}.
\] 
Fix $\xi \in P$ such that $h^0(\tilde{C}, \sO_{\tilde{C}}(2 \tilde{C}) \otimes P_\xi) \geq 3$ and
$D \in | \sO_{\tilde{C}}(2 \tilde{C}) \otimes P_\xi|$.
Let $s$ and $t$ be two sections of $\sO_{\tilde{C}}(\tilde{C})$ such that the associated divisors have
disjoint supports, then we have an exact sequence
\[
0 \rightarrow \sO_{\tilde{C}}(D-\tilde{C}) 
\stackrel{(t,-s)}{\longrightarrow}
\sO_{\tilde{C}}(D)^{\oplus 2} 
\stackrel{(s,t)}{\longrightarrow}
\sO_{\tilde{C}}(D+\tilde{C}) 
\rightarrow 0.
\]
This implies
\[
h^0(\tilde{C}, \sO_{\tilde{C}}(D-\tilde{C})) +  h^0(\tilde{C}, \sO_{\tilde{C}}(D+\tilde{C}))
\geq 2 h^0(\tilde{C}, \sO_{\tilde{C}}(D)) = 6,
\]
furthermore by Riemann-Roch
$h^0(\tilde{C}, \sO_{\tilde{C}}(D+\tilde{C})) = h^0(\tilde{C}, \sO_{\tilde{C}}(K_{\tilde C}-D-\tilde{C})) + 5$.
Now $K_{\tilde{C}} - D \equiv \sigma(D)$ and
$h^0(\tilde{C}, \sO_{\tilde{C}}(\sigma(D)-\tilde{C}))
= h^0(\tilde{C}, \sO_{\tilde{C}}(D- \sigma(\tilde{C})))$ imply 
\[
h^0(\tilde{C}, \sO_{\tilde{C}}(D-\tilde{C})) +  h^0(\tilde{C}, \sO_{\tilde{C}}(D- \sigma(\tilde{C}))) \geq 1.
\]
Hence $D \equiv \tilde{C} + D'$ or $D \equiv \sigma(\tilde{C}) + D'$ where $D'$ is an effective divisor 
such that $\pi_* D' \equiv H$. We see as in the first part of the proof that the effective divisors $D'$ 
such that $\pi_* D' \equiv H$ are parametrised by a set that is contained in a translate of $F \cup (-F)$. 
$\square$

\section{Proof of theorem \ref{theoremfanosurfacegvsheaf}.}
\label{sectionproof}

Since $\sO_F(\Theta) \equiv_{\num} \sO_F(2 \tilde{C})$ by formula (\ref{Theta2C}), 
it is equivalent to verify the stated properties
for the sheaf $\sO_F(2 \tilde{C})$. 

{\it Step 1. The second cohomological support locus is contained in a translate of $F \cup (-F)$.}
By formula (\ref{KF3C}), we have $K_F \equiv \sO_F(3 \tilde{C}) \otimes P_{\xi_0}$ 
for some $\xi_0 \in P$. Hence 
by Serre duality 
$h^2(F, \sO_F(2 \tilde{C}) \otimes P_\xi)=h^0(F, \sO_F(\tilde{C}) \otimes P_\xi^* \otimes P_{\xi_0})$, 
so it is equivalent to consider the non-vanishing locus 
\[
V_0 := \{
\xi \in P \ | \ h^0(F, \sO_F(\tilde{C}) \otimes P_\xi) > 0
\}.
\] 
If $l \in F$ is a line on $X$, the corresponding incidence curve $D_l \subset F$
is an effective divisor numerically equivalent to $\tilde{C}$, so it is clear that
 $\pm F$ is (up to translation) a subset of $V_0$. In order to show that we have an equality, 
consider the exact sequence
\[
0 \rightarrow \sO_F \otimes P_\xi \rightarrow \sO_F(\tilde{C}) \otimes P_\xi
\rightarrow \sO_{\tilde{C}}(\tilde{C}) \otimes P_\xi
\rightarrow 0.
\]
Clearly $h^0(F, \sO_F \otimes P_\xi)=0$ for $\xi \neq 0$, so
$h^0(F, \sO_F(\tilde{C}) \otimes P_\xi) \leq 
h^0(\tilde{C}, \sO_{\tilde{C}}(\tilde{C}) \otimes P_\xi)$ 
for $\xi \neq 0$.
Since a divisor $D \in \ |\sO_{\tilde{C}}(\tilde{C})|$  
satisfies $\pi_* D \equiv H$, we conclude with Lemma \ref{lemmadegreefiveten}.

{\it Step 2. The first cohomological support locus is 
is contained in a union of translate of $F \cup (-F)$.}
Since $\chi(F, \sO_F(2 \tilde{C}))=\chi(F, \sO_F(\Theta))=1$ (formula (\ref{chitheta})),  
we have 
\[
h^1(F, \sO_F(2 \tilde{C}) \otimes P_\xi) = 
h^0(F, \sO_F(2 \tilde{C}) \otimes P_\xi)
+ h^0(F, \sO_F(\tilde{C}) \otimes P_\xi^* \otimes P_{\xi_0})
- 1. 
\]
Since 
\[
h^0(F, \sO_F(2 \tilde{C}) \otimes P_\xi)=h^0(F, \sO_F(\Theta) \otimes P_\xi) \geq 1
\]
for all $\xi \in P$,
the first cohomological support locus is contained in the locus
where $h^0(F, \sO_F(\tilde{C}) \otimes P_\xi^* \otimes P_{\xi_0})>0$ 
or $h^0(F, \sO_F(2 \tilde{C}) \otimes P_\xi)>1$.  
By step 1 the statement follows if we show the following claim:
the set 
\[
V_1 := \{
\xi \in P \ | \ h^0(F, \sO_F(2 \tilde{C}) \otimes P_\xi) > 1
\}
\] 
is contained in a union of translates of $F \cup (-F)$.

{\it Step 3. Proof of the claim and conclusion.}
Consider the exact sequence
\[
0 \rightarrow \sO_F(\tilde{C}) \otimes P_\xi \rightarrow \sO_F(2 \tilde{C}) \otimes P_\xi
\rightarrow \sO_{\tilde{C}}(2 \tilde{C}) \otimes P_\xi
\rightarrow 0.
\]
By the first step we know that
$h^0(F, \sO_F(\tilde{C}) \otimes P_\xi)=0$ for $\xi$ in the complement of a translate of $F \cup (-F)$, 
so
\[
h^0(F, \sO_F(2 \tilde{C}) \otimes P_\xi) \leq 
h^0(\tilde{C}, \sO_{\tilde{C}}(2 \tilde{C}) \otimes P_\xi)
\] 
for $\xi$ in the complement of a translate of $F \cup (-F)$.
The claim is then immediate from Lemma \ref{lemmadegreefiveten}.
By the same lemma $h^0(\tilde{C}, \sO_{\tilde{C}}(2 \tilde{C}) \otimes P_\xi)=1$ for $\xi \in P$ general,
so $h^0(F, \sO_F(2 \tilde{C}) \otimes P_\xi)=h^0(F, \sO_F(\Theta) \otimes P_\xi)=1$
for $\xi \in P$ general.
 $\square$
 
\medskip
 
{\bf Remark.} It is possible to strengthen {\em a posteriori} the statements in the proof:
since Theorem \ref{theoremfanosurfacegvsheaf} holds, we can use the Fourier-Mukai techniques
from \cite{PP06} to see that the cohomological support loci are supported exactly
on the {\it theta-dual} of $F$ (ibid, Definition 4.2), which in our case is just $F$.  

\medskip

{\bf Acknowledgements.} I would like to thank Mihnea Popa for suggesting to me to work on this question. 
Olivier Debarre has shown much patience at explaining to me the geometry of abelian varieties. 
For this and many discussions on minimal cohomology classes I would like to express my deep gratitude.


\begin{thebibliography}{10}

\bibitem{Bea77}
A. Beauville.
\newblock Vari\'et\'es de {P}rym et jacobiennes interm\'ediaires.
\newblock {\em Ann. Sci. \'Ecole Norm. Sup. (4)}, 10(3):309--391, 1977.

\bibitem{Bea82b}
A. Beauville.
\newblock Les singularit\'es du diviseur {$\Theta $} de la jacobienne
  interm\'ediaire de l'hypersurface cubique dans {${\bf P}\sp{4}$}.
\newblock In {\em  Lect. Notes Math.} 947., pages 190--208. Springer, Berlin, 1982.

\bibitem{Bea82}
A. Beauville.
\newblock Sous-vari\'et\'es sp\'eciales des vari\'et\'es de {P}rym.
\newblock {\em Comp. Math.}, 45(3):357--383, 1982.

\bibitem{CG72}
C.~H. Clemens and P.~A. Griffiths.
\newblock The intermediate {J}acobian of the cubic threefold.
\newblock {\em Ann. of Math. (2)}, 95:281--356, 1972.

\bibitem{Deb95}
O. Debarre.
\newblock Minimal cohomology classes and {J}acobians.
\newblock {\em J. Alg. Geom.}, 4(2):321--335, 1995.

\bibitem{H07}
A. H{\"o}ring.
\newblock Paper in preparation.
\newblock {\em Soon on this server}, 2007.

\bibitem{Mum71}
D. Mumford.
\newblock Theta characteristics of an algebraic curve.
\newblock {\em Ann. Sci. \'Ecole Norm. Sup. (4)},
 4:181--192, 1971.


\bibitem{PP03}
G. Pareschi and M. Popa.
\newblock Regularity on abelian varieties. {I}.
\newblock {\em J. Amer. Math. Soc.}, 16(2):285--302, 2003.

\bibitem{PP06}
G. Pareschi and M. Popa.
\newblock {Generic vanishing and minimal cohomology classes on abelian
  varieties}.
\newblock {\em arXiv:math.AG/0610166}, 2006.

\bibitem{PP06b}
G. Pareschi and M. Popa.
\newblock {GV-sheaves, Fourier-Mukai transform, and Generic Vanishing}.
\newblock {\em arXiv:math.AG/0608127}, 2006.

\bibitem{Ran81}
Z. Ran.
\newblock On subvarieties of abelian varieties.
\newblock {\em  Inventiones Math.}, 62:459--479, 1981.

\bibitem{We81}
G.~E. Welters.
\newblock {\em Abel-{J}acobi isogenies for certain types of {F}ano threefolds},
  volume 141 of {\em Mathematical Centre Tracts}.
\newblock Mathematisch Centrum, Amsterdam, 1981.
\end{thebibliography}
\end{document}